\documentclass[12pt]{article}

\usepackage{amssymb, amsmath, url, amssymb, amsthm,color}
\usepackage[T1]{fontenc}
\usepackage[latin1]{inputenc}

\newcommand{\e}{\varepsilon}

\def \Z{{\bf Z}}
\def \qq{{\qquad}}
\def \E{{\bf E}\, }
\def \R{{\bf R}}
\def \T{{\bf T}}
\def \P{{\bf P}}
\def \s{{\sigma}}

\begin{document}

\title{\bf On series $\sum c_k f(kx)$ and
Khinchin's conjecture}

\author{Istv\'an Berkes$^{\rm 1}$ and Michel Weber$^{\rm 2}$ \\[1ex]}

\date{}
\maketitle

\vskip0.2cm

\begin{abstract}
We prove the optimality of a criterion of Koksma (1953) in
Khinchin's conjecture,
settling a long standing open problem in analysis. Using this
result, we also give a near optimal condition for the a.e.\
convergence of series $\sum_{k=1}^\infty c_k f(kx)$ for $f\in L^2$.
\end{abstract}

\renewcommand{\thefootnote}{}
\footnote{$^{\rm 1)}$Institute of Statistics, Graz University of
Technology, Kopernikusgasse 24/III,  8010 Graz, Austria. Email: {\tt
berkes@tugraz.at} . \ Research supported by Austrian Science Fund (FWF),
Grant  P24302-N18 and OTKA Grant K 81928.}

\footnote{$^{\rm 2)}$ IRMA, 10 rue du G\'en\'eral Zimmer,
67084 Strasbourg Cedex, France. Email:
{\tt michel.weber@math.unistra.fr}}

\section{Introduction}\label{se:intro}

Let $\T=\R/\Z\simeq[0,1)$ denote the circle endowed with Lebesgue
measure, $e(x)=\exp(2i\pi x)$, $e_n(x)= e(nx)$, $n\in {\Z}$. Let
\begin{equation}\label{fcond}
f\in L^2(\T),  \quad \int_{\T} f(x)dx=0, \quad f(x)\sim
\sum_{\ell\in {\Z}} a_\ell e_{\ell}, \quad a_0=0.
\end{equation}
Two closely related classical problems of analysis are the almost
everywhere convergence of series
\begin{equation}\label{sum}
\sum_{k=1}^\infty c_k f(kx)
\end{equation}
and the a.e.\ convergence of averages
\begin{equation}\label{avg}
\frac{1}{N}\sum_{k=1}^N f(kx).
\end{equation}
Khinchin \cite{kh} conjectured (assuming only $f\in L^1(\T)$ in
(\ref{fcond})) that
\begin{equation}\label{khin}
\lim_{N\to\infty} \frac{1}{N}\sum_{k=1}^N f(kx)= 0 \qquad
\text{a.e.}
\end{equation}
This conjecture remained open for nearly 50 years. Koksma \cite{ko1}
proved that (\ref{khin}) holds if the Fourier coefficients of $f$
satisfy
\begin{equation}\label{kok0}
\sum_{k=1}^\infty |a_k|^2 (\log\log k)^3<\infty,
\end{equation}
and in \cite{ko2} he weakened the condition to
\begin{equation}\label{kok}
\sum_{k=1}^\infty |a_k|^2 \sigma_{-1} (k)<\infty,
\end{equation}
where
\begin{equation}\label{sigma_s}
\sigma_s (k)=\sum_{d|k} d^s.
\end{equation}
The function $\sigma_{-1} (k)$ is multiplicative and by Gronwall's
estimate \cite{gr} we have
$$
\limsup_{k\to\infty} \frac{\sigma_{-1} (k)}{\log\log k}=e^\lambda,
$$
where $\lambda$ is Euler's constant. Thus condition (\ref{kok}) is
satisfied if
\begin{equation}\label{kok2}
\sum_{k=1}^\infty |a_k|^2 \log\log k<\infty.
\end{equation}
Note the difference between (\ref{kok}) and (\ref{kok2}): by
a theorem of Wintner (\cite{wi}, p.\ 180) the averages
$\frac{1}{J}\sum_{j=1}^J\s_{-1}(j)$ remain bounded, which easily implies
that for any function $\omega(k)\to \infty$ we have $\sigma_{-1} (k)
\ll \omega(k)$ on a set of $k$'s with asymptotic density 1. Thus
(\ref{kok}) is only slightly stronger than $\sum_{k=1}^\infty
|a_k|^2<\infty$, showing that relation (\ref{khin}) holds for
a very large class of functions $f\in L^2({\bf T})$.
However, Marstrand \cite{ma} showed that there exist functions $f\in L^2({\bf T})$
(and even bounded functions $f$) satisfying (\ref{fcond}) such that
(\ref{khin}) fails, thereby disproving Khinchin's conjecture. In his
seminal paper \cite{bo}, Bourgain used his entropy method to
construct a new counterexample and mentioned, concerning conditions
(\ref{kok0}), (\ref{kok2}), that ``{\it ... A more detailed analysis of the previous
construction shows that Koksma's double logarithmic condition is
essentially best possible.}'' See the remark on p.\ 89 of \cite{bo}.
The purpose of the present paper is to verify Bourgain's claim in a
slightly modified form; we will namely prove the following

\bigskip\noindent
{\bf Theorem 1.} {\it Let $w(n)$ be a nonnegative function  of a
natural  argument, which is sub-multiplicative  and bounded in mean.
Assume that
 $$ w(n) =o(\log\log n).$$
Then there exists  a function $f$ satisfying (\ref{fcond}) with
\begin{equation*}
\sum_{k=1}^\infty |a_k|^2 w (k)<\infty
\end{equation*}
such that (\ref{khin}) is not valid.}

\medskip

A function $w$ is called sub-multiplicative if $w(nm)\le w(n)w(m)$
for all $m, n$ with $(m, n)=1$. Clearly, for any $0<\e<1$ the
function $w(n)=\sigma_{-1} (n)^{1-\e}$  satisfies the assumptions of
Theorem 1 and thus it follows that (\ref{fcond}) and
\begin{equation}\label{eps}
\sum_{k=1}^\infty |a_k|^2 \sigma_{-1} (k)^{1-\e}<\infty
\end{equation}
do not generally imply (\ref{khin}). In other words,  Koksma's
condition  (\ref{kok}) is optimal for the a.e.\ convergence relation
(\ref{khin}). Whether Theorem 1 remains valid without the assumption
of sub-multiplicativity and bounded means of $w$
remains open.

\medskip
Theorem 1 shows that $\sigma_{-1} (k)$ is an optimal Weyl factor  in
the Fourier series of $f$ for the validity of (\ref{khin}), but it
does not mean that in the absence of this Weyl factor relation
(\ref{khin}) is always false. In fact, we will see that   under mild
regularity conditions on the Fourier coefficients of $f$, relation
(\ref{khin}) holds under assuming only $f\in L^2$, i.e.\ $\sum_{k=1}
^\infty |a_k|^2<\infty$. Also, the norming factor in (\ref{khin})
can be substantially diminished, see Corollaries 1-3.

\medskip
The previous results give a fairly satisfactory  picture on the
validity of the convergence relation (\ref{khin}). Results
concerning the convergence of sums (\ref{sum}) are much less
complete. Note that, by the Kronecker lemma, the a.e.\ convergence of
(\ref{sum}) with $c_k=1/k$ implies (\ref{khin}), so the two problems
are closely connected. By Carleson's theorem \cite{ca}, in the case
$f(x)=\sin 2\pi x$, $f(x)=\cos 2\pi x$ the series (\ref{sum})
converges a.e.\ provided $\sum_{k=1}^\infty c_k^2<\infty$. Gaposhkin
\cite{ga68} showed that this remains valid if the Fourier series of
$f$ converges absolutely; in particular this holds if $f$ belongs to
the Lip ($\alpha$) class for some $\alpha>1/2$. However, Nikishin
\cite{ni} showed that the analogue of Carleson's theorem fails for
$f(x)=\hbox{sgn} \sin 2\pi x$ and fails also for some continuous
$f$. Recently, Aistleitner and Seip  \cite{aibese} and Berkes
\cite{be2} showed that for $f\in \text{BV}$ the series
(\ref{sum}) converges a.e. provided $\sum_{k=1}^\infty c_k^2
(\log\log k)^\gamma<\infty$ for $\gamma>4$, but generally not for
$\gamma<2$. A similar, slightly weaker result holds for the Lip (1/2)
class, showing, in particular, that Gaposhkin's result above
is sharp. For the class Lip ($\alpha$), $0<\alpha< 1/2$ and other classical
function classes like $C(0, 1)$, $L^p(0, 1)$ or function classes
defined by the order of magnitude of their Fourier coefficients, the
results are much less complete: there are several sufficient
criteria (see e.g. \cite{ai}, \cite{bewe}, \cite{bewe2}, \cite{br},
\cite{ga66a}, \cite{we}; see also \cite{ga66b} for the history of
the subject until 1966) and necessary criteria (see \cite{be}), but
there are large gaps between the sufficient and necessary
conditions. For $f \in \text{Lip} \,(\alpha)$, $0<\alpha<1/2$ Weber
\cite{we} proved that the series (\ref{sum}) converges a.e.\
provided
\begin{equation}\label{dn}
\sum_{k=1}^\infty c_k^2 d(k) (\log k)^2<\infty,
\end{equation}
where $d(n)=\sum_{d|n} 1$ is the divisor function. For functions $f$
with Fourier coefficients of order $O(k^{-s})$, $1/2<s\le 1$ Berkes
and Weber \cite{bewe2} obtained the convergence criterion
\begin{equation}\label{dn2}
\sum_{k=1}^\infty c_k^2 \sigma_{1-2s} (k) (\log k)^2<\infty.
\end{equation}
In view of the role played by the function $\sigma_{-1} (k)$ in Khinchin's
conjecture, the appearance of the arithmetic functions $d(k)$,
$\sigma_{1-2s} (k)$ in (\ref{dn}), (\ref{dn2}) is not surprising.
However, unlike relation (\ref{khin}), there are no comparable
necessary conditions for the convergence of (\ref{sum})
for these function classes: for example, in the case Lip
($\alpha$), $0<\alpha<1/2$  we know only that
\begin{equation}\label{lipalso}
\sum_{k=1}^\infty c_k^2 (\log k)^\gamma<\infty, \qquad 0<\gamma<1-2\alpha
\end{equation}
is in general not sufficient for the a.e.\ convergence of
(\ref{sum})  (see Berkes \cite{be}). Note the large gap between
(\ref{dn}) and (\ref{lipalso}): while the average order of magnitude
of $d(n)$ is $\log n$, we have $d(n)\ge \exp (c\log n/\log\log n)$
for infinitely many $n$. The purpose of this paper is to give a
sufficient condition for the a.e.\ convergence of (\ref{sum}) for
$f\in L^2$ which is optimal up to a logarithmic factor.
Specifically, we will prove the following result.

\bigskip\noindent
{\bf Theorem 2.} {\it Let $f$ satisfy (\ref{fcond}) and put
\begin{equation}\label{conds}
g(r)= \sum_{k=1}^\infty |a_{rk}|^2, \qquad G(r)=\sum_{j\le 2r} g(j),
\qquad h(n)=\sum_{d|n} (dg(d)+G(d)).
\end{equation}
Then $\sum_{k=1}^\infty c_kf(kx)$ converges a.e.\ provided
\begin{equation}\label{ccond}
\sum_{k=1}^\infty c_k^2 h(k) (\log k)^2<\infty.
\end{equation}
On the other hand, for any $\delta>0$ there exists
an $f\in L^2$ and coefficients $c_k$ such that
\begin{equation}\label{nec}
\sum_{k=1}^\infty c_k^2 h(k) (\log k)^{-\delta} <\infty
\end{equation}
but $\sum_{k=1}^\infty c_kf(kx)$ does not converge a.e.}

\bigskip

Theorem 2 is the series analogue of Koksma's theorem, providing a nearly optimal Weyl factor for the a.e.\ convergence of $\sum_{k=1}^\infty c_k f(kx)$ for general $f\in L^2({\bf T})$. The theorem shows that the convergence of $\sum_{k=1}^\infty c_kf(kx)$ is intimately connected with the behavior of the function $g(d)=\sum_{k=1}^\infty |a_{dk}|^2$ which, in turn, depends on the distribution of the
numbers $|a_k|$ on $[0, \infty)$. Under mild regularity
conditions on the $|a_k|$,  $g(d)$ will be
small for large $d$, leading to a small Weyl factor $h$ and better
convergence properties. See Corollaries 1-3 below, where the Weyl factor $h(n)$ reduces
to classical arithmetic functions like $d(n)$, $\sigma_s (n)$. For general
$f\in L^2({\bf T})$, $g(d)$ can behave rather irregularly which, combined
with the summation $d|n$ in the definition of $h(n)$, leads to very irregularly
changing functions $h(n)$. Clearly, $h(n)$ will be small for numbers $n$ with
few prime factors, showing, e.g., that the a.e.\ convergence properties
of dilated sums $\sum_p c_p f(px)$ extended for primes are considerably better
that those of general series (\ref{sum}).

By the Kronecker lemma, Theorem 2 implies for any $f\in L^2({\bf T})$ that
\begin{equation}\label{h2}
\left|\sum_{k=1}^N f(kx)\right| \ll \sqrt{N} (\log N)^{3/2+\varepsilon}
\hat{h} (N)^{1/2} \qquad \text{a.e.}
\end{equation}
where $\hat{h} (n)=\max_{1\le k\le n} h(k)$ is the smallest monotone majorant of the function $h$. 
Using the Kronecker lemma leads to some loss of accuracy, but the corollaries below show that under mild monotonicity or regularity conditions on the $|a_k|$, the right hand side of (\ref{h2}) will be $O(N^{1/2+\e})$, and thus the Khinchin conjecture is valid under such conditions.

\bigskip\noindent
{\bf Corollary 1.} {\it Let $f$ satisfy (\ref{fcond}) where $(a_k)$
satisfies one of the following conditions:

\medskip\noindent
(a) \ $|a_k|$ is regularly varying as $k\to\infty$

\smallskip\noindent
(b) \ $k^{-\gamma} |a_k| $ is non-increasing for some $\gamma>0$

\smallskip\noindent
(c) There exists a $C>0$ such that for any integer $d\ge 1$ we have
$\sum_{k=1}^\infty |a_{dk}|^2\le C/d.$

\medskip\noindent
Then the series $\sum _{k=1}^\infty c_k f(kx)$ converges a.e.\
provided $\sum_{k=1}^\infty c_k^2 d(k) (\log k)^2<\infty$ and
consequently we have
\begin{equation*}\label{st1}
\left|\sum_{k=1}^N f(kx)\right| \ll \sqrt{N} (\log N)^{3/2+\varepsilon}
\hat{d} (N)^{1/2} \qquad \text{a.e.}
\end{equation*}
for any $\varepsilon>0$, where $\hat{d} (N)=\max_{1\le k\le N} d(k)$.}

\bigskip
Since $d(k)\ll k^\e$ for any $\e>0$, we get

\bigskip\noindent
{\bf Corollary 2.} {\it Let $f$ satisfy (\ref{fcond}). Under any of
the regularity conditions (a), (b), (c) we have
$$\left|\sum_{k=1}^N f(kx)\right|\ll N^{1/2+\e} \qquad \text{a.e.}$$
for any $\e>0$.}

\bigskip
It is easily seen that in Corollary 1 both (a) and (b) imply (c), so
(c) is the weakest condition of the three. It is also quite natural:
it requires that all subsums $\sum_{k=1}^\infty |a_{dk}|^2$ carry at
most their ``fair share'' in the sum $\sum_{k=1}^\infty |a_k|^2$.
However, in a number of important cases the estimate $\sum_{k=1}^\infty |a_{dk}|^2
\ll d^{-1}$ can be
improved, leading to a convergence theorem for $\sum_{k=1}^\infty c_k f(kx)$
with a smaller Weyl factor.

\bigskip\noindent
{\bf Corollary 3.} Let $f$ satisfy (\ref{fcond}), where
\begin{equation*} \label{acond}
|a_k|\ll k^{-1/2} \varphi (k) \qquad k=1, 2, \ldots
\end{equation*}
with a non-increasing function $\varphi$ satisfying $\sum_{k=1}^\infty
k^{-1} \varphi^2(k)<\infty$. Let
\begin{equation*}\label{gcond}
\psi(r)=\sum_{k\geq r} k^{-1}\varphi^2(k), \qquad h(N)=\sum_{d|N} \psi
(d)
\end{equation*}
and assume that $\psi$ is regularly varying with exponent $>-1$.
Then the series $\sum_{k=1}^\infty c_k f(kx)$ converges a.e.\ provided
$\sum_{k=1}^\infty c_k^2 h(k) (\log k)^2<\infty$ and consequently (\ref{h2}) holds.

\bigskip
For example, for $g(k)=k^{-\gamma}$, $0<\gamma \le 1/2$
we get $h(k) \ll \sigma_{-2\gamma} (k)$, leading to the convergence condition
(\ref{dn2}) in Berkes and Weber \cite{bewe2}. For
$g(k)=(\log k)^{-\gamma}$, $\gamma >1/2$ we get
$$h(k)\ll \sum_{d|k} (\log d)^{-(2\gamma-1)}.$$

\section{Proof of Theorem 1}

Let $w=\{w_n, n\in \Z\}$ be a sequence of positive reals. Let $L^2_w
$ be the associated Sobolev space on the circle, namely    the
subspace of $L^2  $ consisting with functions $f $ such that
$$ \|f\|_w^2:=\sum_{ n\in \Z} w_n a^2_n(f)<\infty.  $$
This is  a Hilbert space with scalar product defined by  $\langle
f,g\rangle= \sum_{ n\in \Z} w_n a_n(f)a_n(g)$, $f,g\in L^2_w $.

The proof of Theorem 1 is  based on an adaptation to the  Sobolev
space $L^2_w $ of the method elaborated by Bourgain in \cite{bo}.
    Let
 $
f(x)= \sum_{\ell\in {\Z}} a_\ell e_{\ell }   $,  $  a_0=0$ and
consider the    dilation operators $T_jf(x) = f(jx)$.  These are
positive isometries on $L^p $,  $p\ge 1$,     such that $T_j1=1$ for
all $j$,  and for all $f\in L^2 $
$$ \frac{1}{J}\sum_{j\le J} T_j f \buildrel{L^2 }\over {\longrightarrow}
\int f d\l, \qq J\to \infty, $$
To   $f\in L^2 $ we associate
 $$ F_{J,f} =\frac{1}{\sqrt J}\sum_{1\le j\le J} g_j T_jf,\qq (J\ge 1)$$
where $g_1,g_2,\ldots $ are i.i.d. standard Gaussian random variables.

\bigskip\noindent
{\bf Proposition.}  {\it Let   $ S_n\colon L^2 \to L^2 ,\
n=1,2,\dots$ be   continuous operators commuting    with  $T_j$ on
$L^2 $,  $S_nT_j=T_jS_n$  for all $n$ and $j$.
  Assume that the
following property is fulfilled:
$$
\l \Big\{  \sup_{n\ge 1} |S_n(f)| <\infty \Big\}=1,  \qq  \forall  f\in L^2_w    ,
$$
 then   there exists a   constant $C $ depending on   $\{S_n , n\ge 1\}$ only,
such that
$$ \sup_{\e >0} \varepsilon\sqrt{  \log
N_f(\varepsilon)}   \le C \limsup_{J\to \infty}\big(\E \|F_{J, f
}\|_w\big)^{1/2},  \qq  \forall  f\in L^2_w ,$$ where
$N_f(\varepsilon)$ is the entropy number associated with the set
$C_f=\{S_nf, n\ge 1\}$, namely   the minimal   number of $L^2 $ open
balls of radius~$\e$, centered in $C_f$  and enough to cover $C_f$.}

\bigskip\noindent
{\bf Proof.}
 By the Banach principle, there exists a non-increasing function
$C:\R_{+}\to \R_{+}$ such that
$$\forall \varepsilon>0,\forall g\in L_w^2(\T),\qquad \l\Big\{\sup_n|S_n(g)|\ge
\|g\|_{w }C(\varepsilon)\Big\}\le \varepsilon. $$
 Let $0<\e<1/4$. Let $f\in L^2(\T)$.
   Taking $g=F_{J,f}$ and using   Fubini's theorem, gives
$$\int_\T \P\Big\{ \sup_{n\ge 1}|S_n(F_{J,f}) |\ge
C(\varepsilon)\|F_{J,f} \|_{w}\Big \} \ d\l \le \varepsilon.$$
It follows that $$
  \l\left\{ x\in \T: \P\Big\{\omega :\sup_{n\ge
1}|S_n(F_{J,f}(\omega,.))(x)\mid
\ge C(\varepsilon)\|F_{J,f}(\omega,.)\|_{w}\Big\} \ge
\sqrt{\varepsilon}\right\}\le \sqrt{\varepsilon},$$
which is better  rewritten under the following form $$
  \l\left\{x\in \T:\P\Big\{\omega :\sup_{n\ge
1}|S_n(F_{J,f}(\omega,.))(x)\mid
\le C(\varepsilon)\|F_{J,f}(\omega,.)\|_{w}\Big\} \ge
1-\sqrt{\varepsilon}\right\}\ge 1-\sqrt{\varepsilon}.
 $$
By   Tchebycheff's inequality,    $\P\Big\{\| F_{J,f}\|^2_w> \E \| F_{J,f}\|_w^2/  \e\Big\}\le \e  $.
We deduce that
 the set
$$ X_{\varepsilon,J,f}=\Big\{x\in
\T:\P\Big\{\omega:\sup_{n\ge 1}|S_n(F_{J,f}(\omega,.)(x)|\le
  C(\varepsilon)\big(\E \| F_{J,f}\|_w^2/  \e \big)^{1/2}\Big\}\ge
1-2\sqrt{\varepsilon}\Big\} $$
has measure greater than $  1-\sqrt{\varepsilon} $.
 The classical estimate of Gaussian semi-norms implies
$$\forall x\in X_{\varepsilon,J,f},\qquad{ {\bf E}\ }\sup_{n\ge
1}|S_n(F_{J,f}(\omega,))(x)|\le {4
  \over (1-2\sqrt{\varepsilon})}\, {
C(\varepsilon) \over \sqrt{\varepsilon} }\Big(\E \| F_{J,f}\|_w^2  \Big)^{1/2}. $$

Now let $I$ be a finite set of integers such that
$\|S_n(f)-S_m(f)\|_{2 }\not= 0 $,  for all distinct elements $m,n\in
I$. By the commutation property, $S_n(F_{J,f})=F_{J,S_nf} $; so that
\begin{eqnarray*}\E \big|S_n(F_{J,f})-S_m(F_{J,f}\big|^2&=&\E
\big| F_{J,S_nf-S_mf} \big|^2= \frac{1}{J}\sum_{j\le J}( T_j(S_nf-S_mf))^2
\cr &=&\frac{1}{J}\sum_{j\le J}  T_j(S_nf-S_mf) ^2\ \rightarrow \ \|S_nf-S_mf\|_2^2,
\end{eqnarray*}
in $L^2 $ as $J$ tends to infinity. We have used the fact that
$(T_jf)^2= T_jf^2$,  if $f\in L^2 $. By proceeding by extraction, we
  can find a partial index ${\mathcal J}$ such that the set
 $$A(I) =  \bigg\{ \forall J\in {\mathcal J}, \
\forall n,m\in I,\ m\not=n,\quad  {\big( \E |
S_n(F_{J,f})-S_m(F_{J,f}) |^2\big)^{1/2} \over \left\|
S_n(f)-S_m(f)\right\|_{2}} \ge \sqrt{1-\e}  \bigg\},$$ has measure
greater that $  1-\sqrt\e $. \vskip 3pt Let  $J\in {\mathcal J}$,
then $\l( A(I)\cap X_{\varepsilon,J,f})  \ge 1-2\sqrt\e >0$, and for
any $x\in A(I)\cap X_{\varepsilon,J,f}$
$$
C(\varepsilon)\Big(\E \| F_{J,f}\|_w^2  \Big)^{1/2}\ge { {\bf E}\ }\sup_{n\ge
1}|S_n(F_{J,f} )(x)| \ge { {\bf E}\ }\sup_{n\in I} S_n(F_{J,f} )(x) \ge $$
$$ \sqrt{1-\e}\ \E \sup_{n\in I} Z(S_n(f))\ge
(1-2\sqrt\e)\ \ \E \sup_{n\in I} Z(S_n(f)).$$
Therefore,
$$\E \sup_{n\in I} Z(S_n(f)) \le C \Big(\limsup_{J\to \infty}\E \| F_{J,f}\|_w^2  \Big)^{1/2}.$$
Sudakov's minoration implies
$$\sup_{\rho >0} \rho \sqrt{  \log
N_f(\rho)}\le C \Big(\limsup_{J\to \infty} \E \| F_{J,f}\|_w^2  \Big)^{1/2}.$$
\vskip 5 pt

\bigskip\noindent
{\bf Proof of Theorem 1.} Let $P_1,P_2,\dots $ denote the sequence
of prime numbers. Fix some positive integer  $s $ and let $d$ be
some other integer such that $2^d\le P_s$. There exists an integer
$T$ such that if
$$A_T= \{ n=P_1^{\alpha_1}\dots  P_s^{\alpha_s} : 2^T\le n<2^{T+1},\
\alpha_i\ge 0,\ i=1,\dots  ,s \},  $$
then  $\sharp (A_{T+d})\le 2\sharp (A_T)$.
  Put
$$ f=f_T=  {1\over \sharp(A_{T })^{{1/ 2}}}
\sum_{n\in A_{T }}e_n  . $$ It follows from Bourgain's proof
\cite{bo} p.\ 88-89,    (or \cite{WeB} p.\ 239-240 for details)
that
$$N \Big( \big(S_{4^i}(f),i\le  \big[{d\over 2}\big]\big), {1\over 8}\Big)\ge  T.
$$
So that
\begin{equation}\label{min}  \sqrt{  \log T} \le C \Big(\limsup_{J\to \infty} \E \| F_{J,f}\|_w^2  \Big)^{1/2}.
\end{equation} Now as
$$F _{J, f}= \frac{1}{  J^{ 1/ 2 } }\sum_{j\le J } g_j \frac{1}{ \#(A_T)^{ 1/ 2 }}\sum_{n\in A_T }e_{nj}=
\frac{1}{ (J\#(A_T))^{ 1/ 2 }}\sum_{\nu\ge 1}e_\nu   \Big(
\sum_{{1\le j\le J\atop j|\nu}\atop   \frac{\nu}{j}\in A_T }
g_j\Big)$$ we have $$\|F _{J, f}\|_w^2=  \frac{1}{  J\#(A_T)
}\sum_{\nu\ge 1}w_\nu  \Big( \sum_{{1\le j\le J\atop j|\nu}\atop
\frac{\nu}{j}\in A_T}  g_j\Big)^2.$$ These sums are finite sums.
Further,
\begin{eqnarray*}\E \|F_{J, f}\|_w^2&= &\frac{1}{J\#(A_T)}
\sum_{\nu=1}^\infty w_\nu\Big( \sum_{{1\le j\le J\atop j|\nu}\atop
\frac{\nu}{j}\in A_T} 1\Big)= \frac{1}{J\#(A_T)} \sum_{j\le J}
\sum_{m\in A_T}   w_{mj} \cr & \le& \Big(   \frac{1}{J} \sum_{j\le
J} w_{ j}\Big) \Big(  \frac{1}{ \#(A_T)}\sum_{m\in A_T}   w_{m }
\Big)\cr & \le&
   \Big(   \frac{1}{J } \sum_{j\le J} w_{ j}\Big)  \max_{m\in A_T} w_m
\end{eqnarray*}
Therefore,
$$ \limsup_{J\to \infty}\E \|F_{J, f}\|_w^2 \le     \Big(   \limsup_{J\to \infty}
\frac{1}{J } \sum_{j\le J} w_{ j}\Big)  \max_{m\in A_T} w_m \le M  \max_{m\in A_T} w_m, $$
where $M<\infty$  and further
$$ \max_{m\in A_T} w_m= o \big( \max_{m\in A_T} \log\log m \big) = o(\log T),$$
by assumption. Consequently
\begin{equation}\label{maj}  \limsup_{J\to \infty}\big(\E \|F_{J, f_T}\|_w\big)^{1/2}= o(\sqrt{\log T}).
\end{equation}
But this contradicts (\ref{min}), completing the proof of Theorem 1.

\bigskip
\section{Proof of Theorem 2.}

Clearly it suffices to prove Theorem 2 for real valued $f$, when for the Fourier coefficients we have
$a_{-\ell}=\overline{a}_\ell$, $\ell \in {\mathbb Z}$. We first prove the following lemma.

\bigskip\noindent
{\bf Lemma.} Let $f$ satisfy (\ref{fcond}). Then for any $r\ge 1$
and any real coefficients $c_j$ we have
\begin{equation}\label{lem}
\int_0^1 \left(\sum_{\ell=2^r+1}^{2^{r+1}} c_\ell f(\ell
x)\right)^2 dx  \le \sum_{\ell=2^r+1}^{2^{r+1}} c_\ell^2 h(\ell),
\end{equation}
where the arithmetic function $h$ is defined by (\ref{conds}).

\bigskip\noindent
{\bf Proof.} Fix $m, n \ge 1$ and put $m'=m/d, \, n'=n/d$, where
$d=(m, n)$.
Using (\ref{fcond}) and $a_{-\ell}=\overline{a}_\ell$, we get
\begin{align}\label{cov}
&\lambda_{m, n}:=\left|\int_0^1
f(mx)f(nx)dx\right|
 =\left|\sum_{mk=nl \atop k, l\in {\mathbb Z}} a_k\overline{a}_l \right|
 \le \sum_{mk=nl \atop k, l\in {\mathbb Z}} |a_k||a_l|\nonumber \\
 &=2\sum_{mk=nl \atop k, l\ge 1} |a_k||a_l| =2\sum_{m'k=n'l \atop k, l\ge 1} |a_k||a_l|.
\end{align}
Since $(m', n')=1$, the equation $m'k=n'l$ implies that $m'$ is a
divisor of $l$, i.e.\ $l=m'i$ and consequently $k=n'i$ for some
$i\ge 1$. Thus the last expression  in (\ref{cov}) equals
\begin{align} \label{m'n'}
&2\sum_{i=1}^\infty |a_{m'i}||a_{n'i}|\le
\sum_{i=1}^\infty (|a_{m'i}|^2 + |a_{n'i}|^2)=g(m')+ g(n').
\end{align}
Now for any $r\ge 1$ and any coefficients $c_\ell$,
\begin{align}\label{MN}
&\int_0^1 \left(\sum_{\ell=2^r+1}^{2^{r+1}} c_\ell f(\ell
x)\right)^2 dx \nonumber \\
&\le \sum_{i, j=2^r+1}^{2^{r+1}} \lambda_{i, j} |c_i||c_j|\le
\frac{1}{2} \sum_{i, j=2^r+1}^{2^{r+1}} \lambda_{i, j} (c_i^2+c_j^2)
= \sum_{i, j=2^r+1}^{2^{r+1}} \lambda_{i, j} c_i^2
= \sum_{i=2^r+1}^{2^{r+1}} c_i^2 \rho (i)
\end{align}
where
\begin{equation}\label{rho}
\rho (i)= \sum_{j=2^k+1}^{2^{k+1}} \lambda_{i, j}\qquad \text{for}
\ 2^k< i\le 2^{k+1}.
\end{equation}
Thus using (\ref{cov}), (\ref{m'n'}) we get for $2^k< i\le 2^{k+1}$,
\begin{align*}\label{rhoi}
\rho (i)= \sum_{j=2^k+1}^{2^{k+1}} \lambda_{i, j}\le
\sum_{j=2^k+1}^{2^{k+1}}(g (i/(i, j))+ g( j/(i, j)).
\end{align*}
Fix $i$ and $d|i$ and sum here for all $j$ with $(i, j)=i/d$. Then
$j=ri/d$ for some $r\le 2d$ and thus the contribution of these
terms is
\begin{equation*}
\ll \sum_{r\le 2d} (g(d) + g(r))\ll dg(d) +G(d).
\end{equation*}
Thus summing now for $d|i$, we get
$$\rho(i)\ll \sum_{d|i} (dg(d)+G(d))=h(i). $$
The lemma now follows from (\ref{MN}).

\bigskip\noindent{\bf Proof of Theorem 2.} Using the Lemma, the proof
of the sufficiency part can be completed by using the method of
Rademacher and Mensov, see e.g. \cite{al}, pp.\ 80--81.
By the Lemma and (\ref{ccond}) we have
\begin{align*}
&\sum_{r=1}^\infty \int_0^1 r^2\left[ \sum_{j=2^r+1}^{2^{r+1}}c_j
f(jx)\right]^2 dx\ll
\sum_{r=1}^\infty r^2 \sum_{j=2^r+1}^{2^{r+1}}c_j^2 h (j)\\
& \ll \sum_{r=1}^\infty \sum_{j=2^r+1}^{2^{r+1}} c_j^2 (\log j)^2
h (j)<\infty.
\end{align*}
Thus
$$
\sum_{r=1}^\infty r^2\left[ \sum_{j=2^r+1}^{2^{r+1}}c_j
f(jx)\right]^2< \infty \qquad \text{a.e.}
$$
and the Cauchy-Schwarz inequality yields for any $1\le M<N$
\begin{align*}
& \left|\sum_{j=2^M+1}^{2^N} c_j f(jx)\right|^2 \le \left(
\sum_{k=M}^{N-1} \left|\sum_{j=2^k+1}^{2^{k+1}} c_j
f(jx)\right|\right)^2 \nonumber \\
&\le \left(\sum_{k=M}^{N-1} \frac{1}{k^2}\right) \left(
\sum_{k=M}^{N-1} k^{2} \left|\sum_{j=2^k+1}^{2^{k+1}} c_j
f(jx)\right|^2 \right)
 \le 2\sum_{k=M}^\infty k^{2} \left|\sum_{j=2^k+1}^{2^{k+1}} c_j
f(jx)\right|^2 \to 0
\end{align*}
as $M\to\infty$. This implies that $\sum_{j=1}^{2^m} c_j f(jx)$
converges a.e.\ as $m\to\infty$. Now the Lemma and standard
maximal inequalities (see e.g.\  \cite{WeB}, Lemma 8.3.4) imply that
\begin{align*}
&\sum_{k=1}^\infty \left\| \max_{2^k+1\le i\le j\le 2^{k+1}}
|\sum_{\ell=i}^j c_\ell f(\ell x)|\right\|^2  \ll \sum_{k=1}^\infty
k^2 \left(\sum_{\ell=2^k+1}^{2^{k+1}} c_\ell^2 h (\ell)\right)\ll
\sum_{\ell=1}^\infty c_\ell^2 (\log \ell)^2 h(\ell) <\infty
\end{align*}
which yields
\begin{equation}\label{e}
\max_{2^k+1\le i\le j\le 2^{k+1}} \left|\sum_{\ell=i}^j c_\ell
f(\ell x)\right|\to 0 \qquad \text{a.e.}
\end{equation}
proving the first part of Theorem 2.

To prove the second statement of Theorem 2, note that for
$c_k=1/k$ and for any positive non-increasing, slowly varying sequence $(\e_k)$ we
have $\sum_{k=1}^\infty c_k^2 h(k) \e_k= \sum_1+\sum_2$, where
\begin{align}\label{sigma1}
&{\sum}_1= \sum_{k=1}^\infty k^{-2} \e_k \sum_{d|k}
dg(d)=\sum_{d=1}^\infty dg(d)\sum_{j=1}^\infty (dj)^{-2} \e_{dj}\ll \sum_{d=1}^\infty
dg(d) d^{-2} \e_d \nonumber \\
&=\sum_{d=1}^\infty
\frac{\e_d}{d} \sum_{k=1}^\infty |a_{dk}|^2=\sum_{j=1}^\infty
|a_j|^2\sum_{d|j} \frac{ \e_d}{d} =\sum_{j=1}^\infty |a_j|^2
\widetilde{\sigma} (j),
\end{align}
with
\begin{equation}\label{kok3}
\widetilde{\sigma} (k)=\sum_{d|k} \e_d /d.
\end{equation}
Similarly,
\begin{align}\label{sigma2}
&{\sum}_2= \sum_{k=1}^\infty k^{-2} \e_k \sum_{d|k}
G(d)=\sum_{d=1}^\infty G(d)\sum_{j=1}^\infty (dj)^{-2} \e_{dj}\ll \sum_{d=1}^\infty
G(d) d^{-2} \e_d \nonumber \\
&\ll \sum_{d=1}^\infty g(d) d^{-1} \e_d,
\end{align}
which is the same bound as the last expression in the first line of (\ref{sigma1}) and thus
continuing, we get the same estimate as in (\ref{sigma1}). To justify the last step in (\ref{sigma2}),
set $G(0)=0$, $S_d=\sum_{j=d}^\infty \e_j j^{-2}$ and note that
\begin{align*}
&\sum_{d=1}^\infty G(d) d^{-2} \e_d = \sum_{d=1}^\infty G(d)(S_d-S_{d+1})= \sum_{\ell=1}^\infty (G(\ell)-G(\ell-1))S_\ell \\
&=\sum_{\ell=1}^\infty (g(2\ell)+g(2\ell-1)) S_\ell \ll  \sum_{\ell=1}^\infty (g(2\ell)+g(2\ell-1)) \frac{\e_\ell}{\ell}\\
& \ll \sum_{\ell=1}^\infty g(2\ell) \frac{\e_{2\ell}}{2\ell} + \sum_{\ell=1}^\infty g(2\ell-1) \frac{\e_{2\ell-1}}{2\ell-1}
=\sum_{r=1}^\infty g(r) \frac{\e_r}{r}.
\end{align*}
Thus we proved
\begin{equation}\label{twosums}
\sum_{k=1}^\infty c_k^2 h(k)\e_k \ll \sum_{j=1}^\infty a_j^2 \widetilde{\sigma} (j).
\end{equation}
Now choosing $\varepsilon_k=(\log k)^{-\delta}$ we have
\begin{align}\label{sigma}
&\widetilde{\sigma}(k)=\sum_{d|k} \frac{\e_d}{d}=\sum_{d|k, \, d\le \exp( \sigma_{-1} (k)^\delta)} \frac{\e_d}{d}
+\sum_{d|k, \, d> \exp( \sigma_{-1} (k)^\delta)} \frac{\e_d}{d}\ll \sum_{d\le \exp( \sigma_{-1} (k)^\delta)}
\frac{1}{d}   \nonumber\\
&+
\varepsilon_{\exp( \sigma_{-1} (k)^\delta)}\sum_{d|k}\frac{1}{d}
\ll  (\sigma_{-1} (k))^\delta + (\sigma_{-1} (k))^{-\delta^2} \sigma_{-1} (k)
\ll (\sigma_{-1} (k))^{1-\delta^2}.
\end{align}
By Theorem 1 we can choose a function $f$ satisfying (\ref{fcond}) such that
$\sum_{j=1}^\infty |a_j|^2 \sigma_{-1} (j)^{1-\delta^2}$ converges,
but $N^{-1}\sum_{k=1}^N f(kx)$ does not converge a.e. But then by relations (\ref{twosums}) and (\ref{sigma})
we have $\sum_{k=1}^\infty c_k^2 h(k)\e_k <\infty$ for $c_k=1/k$ and $\sum_{k=1}^\infty c_k f(kx)$ cannot
converge a.e., since then by the Kronecker lemma we would have $N^{-1}\sum_{k=1}^N f(kx)\to 0$ a.e.

\bigskip
\bigskip\noindent {\bf Proof of Corollary 3.} Assuming the regular variation of $g(n)$, the statement is an easy
consequence of Theorem 2. However, specializing the proof of Theorem 2 to this case instead, we get the statement without any additional assumptions.

Extend $\varphi$ to $[1, +\infty)$ in a monotone non-increasing fashion.
Then the contribution of the terms in the first sum in (\ref{m'n'}) for $i\ge 2$ can be estimated for $m\le n$
as follows:
\begin{align}\label{lambda}
&2\sum_{i=2}^\infty |a_{m'i}||a_{n'i}|\ll
(m'n')^{-1/2}\sum_{i=2}^\infty i^{-1}\varphi (m'i) \varphi (n'i) \ll
(m'n')^{-1/2}\sum_{i=2}^\infty i^{-1}\varphi ^2(m'i) \nonumber \\
&\ll (m'n')^{-1/2}\int_1^\infty x^{-1} \varphi^2(m'x)\, dx=
(m'n')^{-1/2}\int_{m'}^\infty y^{-1} \varphi^2(y)\, dy \nonumber \\
&\ll (m'n')^{-1/2}\sum_{k\ge m'} k^{-1} \varphi^2(k).
\end{align}
On the other hand, well known properties of regularly varying functions
(see e.g.\ \cite{bgt}, Theorem 1.5.11 (ii) with $f(x)=\varphi^2(x)$, $\sigma=-1$)
imply that the ratios
\begin{equation}\label{ratios}
\varphi^2 (r) /\int_r^\infty t^{-1}\varphi^2 (t) dt \quad \text{and} \quad
\varphi^2 (r) /\sum_{k=r}^\infty k^{-1}\varphi^2 (k)
\end{equation}
converge, as $r\to\infty$, to
a finite limit $c\ge 0$.
Thus for $m\le n$ we have
$$  |a_{m'}||a_{n'}| \ll (m'n')^{-1/2}\varphi(m')\varphi (n') \ll (m'n')^{-1/2}\varphi^2(m')
\ll (m'n')^{-1/2} \sum_{k=m'}^\infty k^{-1}\varphi^2 (k).
$$
Hence (\ref{lambda}) implies for $m\le n$, adding the term for $i=1$,
\begin{align}\label{lambda2}
&2\sum_{i=1}^\infty |a_{m'i}||a_{n'i}|
\ll (m'n')^{-1/2}\sum_{k\ge m'} k^{-1} \varphi^2(k)
= (m'n')^{-1/2} \psi (m') \nonumber\\
&=(m'n')^{-1/2} \psi (m'\wedge
n') =(mn)^{-1/2} (m, n) \, \psi( (m\wedge n)/(m, n))
\end{align}
and the same estimate holds for $m\ge n$. Recall now that for any $r\ge 1$
and any coefficients $c_\ell$
we have (\ref{MN}), where $\rho (i)$ is defined by (\ref{rho}).
Since $\lambda_{m, n}\ll \sum_{i=1}^\infty |a_{m'i}||a_{n'i}|$, using
(\ref{lambda2}) we get for $2^k< i\le 2^{k+1}$, using the
monotonicity and regular variation of $\psi$,
\begin{align*}
\rho (i)= \sum_{j=2^k+1}^{2^{k+1}} \lambda_{i, j}\ll
\sum_{j=2^k+1}^{2^{k+1}}(ij)^{-1/2} (i, j)\psi( (i\wedge j)/(i,
j))\ll 2^{-k} \sum_{j=2^k+1}^{2^{k+1}}(i, j)\psi(j/(i, j)).
\end{align*}
Fix $i$ and $d|i$ and sum here for all $j$ with $(i, j)=i/d$. Then
$j=ri/d$ for some $r\le 2d$ and thus the contribution of these
terms is
\begin{equation*}
\ll 2^{-k} \sum_{r\le 2d} (i/d) \psi\left( \frac{ri/d}{i/d}\right)
\ll (1/d) \sum_{r\le 2d} \psi(r)\ll \psi(d),
\end{equation*}
where we used the regular variation of $\psi$ with index $>-1$ and Theorem 1.5.11 of
\cite{bgt},  p.\ 28 with $\sigma=0$. Thus summing now for $d|i$, we get
$$\rho(i)\ll \sum_{d|i} \psi(d)=h(i), $$
and finally by (\ref{MN})
\begin{align}\label{ut}
\int_0^1 \left(\sum_{\ell=2^r+1}^{2^{r+1}} c_\ell f(\ell
x)\right)^2 dx  \le \sum_{\ell=2^r+1}^{2^{r+1}} c_\ell^2 h(\ell).
\end{align}
The proof can now be completed as in Theorem 2.

\bigskip\noindent
{\bf Proof of Corollary 1.} Assume condition (c) of the Corollary.
Then the first expression in (\ref{m'n'}) can be bounded as
\begin{equation*}\label{biz}
2\sum_{i=1}^\infty |a_{m'i}||a_{n'i}|\le
2\left(\sum_{i=1}^\infty |a_{m'i}|^2\right)^{1/2}
\left(\sum_{i=1}^\infty |a_{n'i}|^2\right)^{1/2}\le C (m'n')^{-1/2}
\end{equation*}
and thus in this case $\lambda_{m, n}$ can be bounded by the last
expression of (\ref{lambda2}) with $\psi = 1$.  Hence the rest of
the proof of Corollary 3 applies with $\psi=1$, showing that (\ref{ut}) holds
with $h(n)=\sum_{d|n}1 =d(n)$. Following the proof of Theorem 2, the statement of Corollary 1
follows in the case (c).

\bigskip
Next we show that in Corollary 1 we have (a)$\Longrightarrow$ (c)
and (b)$\Longrightarrow$ (c). Assume first that (b) holds, then
$$ |a_{n+1}/a_n|\le |(n+1)/n|^\gamma \le 1+ C/n \qquad (n\ge 1)$$
for some constant $C>0$. Let now $k\ge 1, d\ge 2$ and $0\le j \le
d/2$. Then
we get, setting $C_1=e^C$,
\begin{align*}
&|a_{kd}/a_{kd-j}|\\
=&\prod_{r=kd-j}^{kd-1}|a_{r+1}/a_{r}|\le
\prod_{r=kd-j}^{kd-1}(1+C/r)\le \exp \left(\sum_{r=kd-j}^{kd-1}
C/r \right)\le \exp (Cj/(d/2)) \le C_1
\end{align*}
and consequently
$$
\sum_{n=1}^\infty |a_n|^2 \ge \sum_{k=1}^\infty \sum_{j=1}^{[d/2]}
|a_{kd-j}|^2 \ge [d/2] \, C_1^{-2} \sum_{k=1}^\infty |a_{kd}|^2,
$$
proving the
validity of condition (c). If condition (a) holds, then by
$\sum_{k=1}^\infty |a_k|^2<\infty$ its exponent of regularity is negative,
i.e.\ there exists a $\rho>0$ such that $n^\rho |a_n| $ is slowly varying.
But then by the remark in
\cite{bgt}, p.\ 23 preceding Theorem 1.5.4, there exists a non-increasing
sequence $b_n\sim |a_n|$. Clearly, $\sum_{n=1}^\infty b_n^2<\infty$ and by
monotonicity, $(b_n)$ satisfies condition (c). But then $(a_n)$ also
satisfies condition (c).

\end{document}